\documentclass[reqno,a4paper,12pt]{amsart} 

\usepackage{amsmath,amscd,amsfonts,amssymb}
\usepackage{mathrsfs,dsfont}

\allowdisplaybreaks

\numberwithin{equation}{section}
\numberwithin{figure}{section}

\newcommand\R{\mathbb{R}}

\newcommand\Z{\mathbb{Z}}

\newcommand\del{\delta}

\newcommand\lam{\lambda}
\newcommand\Lam{\Lambda}

\newcommand\Om{\Omega}

\newcommand\eps{\varepsilon}
\newcommand{\pphi}{\varphi}

\renewcommand\le{\leqslant}
\renewcommand\ge{\geqslant}
\renewcommand\leq{\leqslant}

\newcommand\sbt{\subset}

\newcommand{\ft}[1]{\widehat{#1}}
\newcommand{\dotprod}[2]{\langle #1 , #2 \rangle}
\newcommand{\mes}{\operatorname{mes}}

\newcommand{\supp}{\operatorname{supp}}

\newcommand{\diam}{\operatorname{diam}}
\newcommand{\dist}{\operatorname{dist}}
\newcommand{\sign}{\operatorname{sign}}

\theoremstyle{plain}
\newtheorem{thm}{Theorem}[section]
\newtheorem{lem}[thm]{Lemma}
\newtheorem{lemma}[thm]{Lemma}

\newtheorem*{claim*}{Claim}

\newcommand{\thmref}[1]{Theorem~\ref{#1}}
\newcommand{\secref}[1]{Section~\ref{#1}}

\newcommand{\lemref}[1]{Lemma~\ref{#1}}

\theoremstyle{definition}
\newtheorem{definition}[thm]{Definition}
\newtheorem*{definition*}{Definition}
\newtheorem*{remarks*}{Remarks}
\newtheorem*{remark*}{Remark}

\newenvironment{enumerate-alph}
{\begin{enumerate}
\addtolength{\itemsep}{5pt}
}
{\end{enumerate}}

\newenvironment{enumerate-num}
{\begin{enumerate}
\addtolength{\itemsep}{5pt}
}
{\end{enumerate}}

\newenvironment{enumerate-text}
{\begin{enumerate}
\addtolength{\itemsep}{5pt}
}
{\end{enumerate}}

\begin{document}

\title
[Unconditional Schauder frames in $L^p$ spaces]
{Unconditional Schauder frames of exponentials and of uniformly bounded functions in $L^p$ spaces}

\author{Nir Lev}
\address{Department of Mathematics, Bar-Ilan University, Ramat-Gan 5290002, Israel}
\email{levnir@math.biu.ac.il}

\author{Anton Tselishchev}
\address{St. Petersburg Department of Steklov Mathematical Institute, Fontanka 27, St. Petersburg 191023, Russia}
\email{celis\_anton@pdmi.ras.ru}

\date{September 9, 2025}
\thanks{N.L.\ is supported by ISF Grant No.\ 1044/21. 
A.T.\ is supported by the Foundation for the 
Advancement of Theoretical Physics and
Mathematics ``BASIS''}

\begin{abstract}
It is known that there is no unconditional basis of exponentials in the space $L^p(\Omega)$, $p \ne 2$, for any set $\Omega \subset \mathbb{R}^d$ of finite measure. This is a consequence of a more general result due to Gaposhkin, who proved that the space $L^p(\Omega)$ does not admit a seminormalized unconditional basis consisting of uniformly bounded functions. We show that the latter result fails if the word ``basis'' is replaced with ``Schauder frame''. On the other hand we prove that if $\Omega$ has nonempty interior then there are no unconditional Schauder frames of exponentials in the space $L^p(\Omega)$, $p \ne 2$.
\end{abstract}

\maketitle


\section{Introduction}

\subsection{}
The most basic fact which lies at the foundation of Fourier analysis is that  the system of exponential functions  $\{e^{2\pi i n x}\}$, $n\in\Z$, forms an orthogonal basis in the Hilbert space $L^2[0,1]$. More generally, for any dimension $d$ the system  $\{e^{2\pi i \dotprod{n}{x}}\}$, $n\in\Z^d$, forms an orthogonal basis in  $L^2$ on the unit cube $[0,1]^d$.

A natural question arises as to which domains $\Omega\subset\R^d$ admit
a countable set of frequencies $\Lambda\subset\R^d$
such that the system of exponential functions
\begin{equation}
\label{eqI1.1}
E(\Lambda)=\{e_\lambda, \; \lambda\in \Lambda\}, \quad e_\lambda(x)=e^{2\pi
i\dotprod{\lambda}{x}},
\end{equation}
forms an orthogonal basis in
 the space $L^2(\Omega)$.
 This question goes back to Fuglede \cite{Fug74} who famously conjectured that there exists an orthogonal basis of exponentials in $L^2(\Om)$ if and only if one can tile $\R^d$ by translated copies of $\Omega$. Fuglede's conjecture inspired extensive research over the years 
and a number of positive as well as negative results have been obtained,
see e.g.\ \cite{IKT01}, \cite{KP02}, \cite{Tao04}.
A major recent result states that the conjecture holds for convex domains in all dimensions 
\cite{LM22}. In particular,  many reasonable domains,
such as  a disk or a triangle in the plane,  do not admit
 an orthogonal basis of exponentials.
We refer the reader to the comprehensive survey \cite{Kol24} for the history of the problem, an overview of the known results and the relevant references.

\subsection{}
If a domain $\Om \sbt \R^d$
 does not have an orthogonal basis of exponentials,
then a \emph{Riesz basis} is the next best thing one can look for.
A system of vectors $\{u_n\}$ in a separable Hilbert space $H$ 
is called a Riesz basis if it can be obtained as the image of an orthonormal basis 
under a bounded and invertible linear operator. Riesz bases 
 share many of the qualities of orthonormal bases.  In particular,
every element $x\in H$ admits a unique series expansion
$x=\sum_{n}c_n u_n$, and moreover the series converges
unconditionally, i.e.\ it converges for any rearrangement of its terms
(see \cite[Section 1.8]{You01}).

The problem as to whether the space $L^2(\Omega)$ admits
a Riesz basis of exponentials $E(\Lam)$ is not well understood, and even for 
 some natural domains $\Om$ such as a disk or a triangle in the plane  the answer
is not known. 	 It is known that  there exists a  Riesz basis  of exponentials for
 any finite union of intervals in $\R$ \cite{KN15}
  (see also \cite{KN16} for a multi\-dimensional version), as well as
for  any zonotope in $\R^d$, namely, any centrally symmetric convex polytope
all of whose faces of all dimensions are also centrally symmetric  \cite{DL22}. 
 An example of a set $\Om$ 
 such that the space $L^2(\Omega)$ does not have a Riesz basis of exponentials was constructed only recently in \cite{KNO23} (in this example, $\Om \sbt \R$ is the union of countably many intervals  with a single accumulation point).

\subsection{}
If there is no Riesz basis of exponentials in
$L^2(\Omega)$ (or we do not know whether one exists)
then we can consider \emph{frames} of exponentials.
A system of vectors $\{u_n\}$  in a separable Hilbert space $H$ 
is called a frame (in the sense of  Duffin and Schaeffer \cite{DS52})
if there exist positive constants $A,B$  such that the inequalities
\begin{equation}\label{eq:frame_ineq}
A \|x\|^2 \le \sum_n |\langle x, u_n \rangle|^2\le B\|x\|^2
\end{equation}
hold for every $x\in H$. In this case there exists another frame
$\{v_n\}\subset H$ such that 
\begin{equation}\label{eq:frame_convergence}
x=\sum_{n} \dotprod{x}{v_n}u_n 
\end{equation}
for every $x \in H$,
and the series converges unconditionally
(see \cite[Section 4.7]{You01}).

We note that if
$\{u_n\}$  is a Riesz basis in $H$ then it is a frame, but in general
a frame need not be a Riesz basis and the
system $\{v_n\}$ satisfying  \eqref{eq:frame_convergence}
need not be unique.

 It was proved in \cite{NOU16} that for any set $\Omega\subset\R$ of finite measure (bounded or unbounded) there exists a set $\Lambda\subset\R$ such that the exponential system $E(\Lambda)$ forms a frame in the space $L^2(\Omega)$.  The same is true also in $\R^d$, see
 \cite[Remark 1]{NOU16}.

\subsection{}
In the present paper we consider similar problems in $L^p$ spaces with $p\neq 2$. 
We start with a few definitions.

A system $\{u_n\}$ in a separable Banach space $X$ is called an \emph{unconditional (Schauder) basis} if every $x\in X$ admits a unique series expansion $x = \sum_n c_nu_n$,
and this series converges unconditionally. The most well-known example of an unconditional basis in the space $L^p[0,1]$, $1 < p < +\infty$, is the Haar basis
(see \cite[II.B.9--10, II.D.13]{Woj91}).

We say that a system $\{u_n\}\subset X$ is \emph{seminormalized}
 if there are positive constants $A,B$  such that 
 $A \le\|u_n\|\le B$ for every $n$. 

It is well-known that  a system $\{u_n\}$  in a Hilbert space $H$ 
is a	 Riesz basis if and only if it is
seminormalized and forms an unconditional basis in $H$.

The classical exponential system with integer frequencies $E(\Z)$
\emph{does not} form an unconditional basis in the space $L^p[0,1]$,
$p \ne 2$, see \cite[Proposition II.D.9]{Woj91}.
We note, however, that the system can be ordered so as to
form a (not unconditional) Schauder basis in  $L^p[0,1]$, $1<p<+\infty$,
see \cite[II.B.11]{Woj91}.

Moreover, for any set
$\Omega\subset\R^d$ of positive and finite measure, 
there does not exist any unconditional basis of exponentials
in the space $L^p(\Omega)$, $p \ne 2$.
 This is a consequence of a more general result due to Gaposhkin  \cite{Gap58},
 who proved that the space $L^p(\Omega)$,  $p \ne 2$, does not admit a  
seminormalized unconditional basis consisting of uniformly bounded functions.
Recall that a system of functions $\{g_n\} \sbt L^p(\Omega)$  is said to be
\emph{uniformly bounded} if there is a constant $M$ such that
for every $n$ we have $|g_n(x)| \le M$ a.e.

For a certain generalization of Gaposhkin's result, see \cite[Theorem 1.2]{AACV19}.

We note that the seminormalization condition in this result is essential: otherwise we could simply take the Haar basis and multiply each element by a sufficiently small positive number, in order to get an unconditional basis in $L^p[0,1]$, $1<p<+\infty$, consisting of uniformly bounded functions.

\subsection{}
There is also a generalization of the notion of a frame 
to the setting of Banach spaces.
If  $X$ is a separable Banach space with dual space $X^*$, then a system of elements
$\{(u_n, u_n^*)\}_{n=1}^{\infty} \sbt X\times X^*$ is called a \emph{Schauder frame} (or a \emph{quasibasis}) if every $x\in X$ admits a series expansion
\begin{equation}
\label{eq:R1.1}
x=\sum_{n=1}^{\infty} u_n^*(x) u_n.
\end{equation}
If the series \eqref{eq:R1.1} 
converges unconditionally for every $x \in X$, then
$\{(u_n, u_n^*)\}_{n=1}^{\infty}$
is called an \emph{unconditional Schauder frame}.

Note that if $\{u_n\}$ is an
unconditional basis in $X$, then there exists a unique system 
of biorthogonal  coefficient functionals
$\{u_n^*\} \sbt X^*$
such that \eqref{eq:R1.1} holds for every $x \in X$
(see \cite[Corollary II.B.7]{Woj91}). 	In this sense, every 
unconditional basis is an unconditional Schauder frame.
However, in general the coefficient functionals
$\{u_n^*\}$ need not be unique
 and they are not necessarily biorthogonal to $\{u_n\}$.

Since for any set
$\Omega\subset\R^d$ of positive and finite measure
there are no unconditional bases of exponentials
in the space $L^p(\Omega)$, $p \ne 2$, 
 the question arises as to whether one can find
at least an unconditional Schauder frame of exponentials
in the space.

The answer would be `no'  if one could extend Gaposhkin's result 
from the context of unconditional bases to unconditional 
Schauder frames, namely, if one could prove
 that the space $L^p(\Omega)$, $p \ne 2$,  does not admit a  
seminormalized unconditional Schauder frame consisting of uniformly bounded functions.

However, our first theorem shows that 
Gaposhkin's result fails if the word ``basis'' is replaced 
with ``Schauder frame''. Moreover, there exists 
an unconditional Schauder frame 
in the space $L^p(\Omega)$, $1<p<+\infty$,  
consisting of unimodular functions:

\begin{thm}
\label{thm:bounded_unc_frames}
Let $\Omega\subset\R^d$ be a set of positive and finite measure.
Then for any $p>1$ there exists an unconditional Schauder frame 
$\{(g_n, g_n^*)\}$ in the space $L^p(\Omega)$ 
such that the functions $g_n$ are unimodular,
i.e.\  satisfy  the condition $|g_n(x)| = 1$ a.e.
\end{thm}

In fact, we will prove two versions of this result. In the first version
the functions $g_n$ are \emph{real-valued} and unimodular,
i.e.\ they are $\{-1, +1\}$-valued. In the second version,
the functions $g_n$ are \emph{smooth} (complex-valued)  and unimodular.

We note that \thmref{thm:bounded_unc_frames}
cannot hold for $p=1$, since in the space $L^1(\Omega)$ 
there are no unconditional Schauder frames, see e.g.\ \cite[Section 4.3]{BC20}.

\subsection{}

\thmref{thm:bounded_unc_frames} suggests that the problem of existence of unconditional Schauder frames of exponentials in $L^p(\Om)$  is more subtle than the corresponding question for unconditional bases: it cannot be settled in the negative by general considerations and   the answer must involve more specific properties of the exponential functions.

 In our second result, we provide an almost complete answer to this question.

\begin{thm}\label{thm:noframes_exp}
Let $\Omega\subset\R^d$ be a set of positive and finite measure. Then,
\begin{enumerate-num}
\item \label{nsf:i}
If $1 < p < 2$ then the space $L^p(\Omega)$ does not admit
any unconditional Schauder frame of exponentials, i.e.\ of the form 
$\{(e_\lambda, e_\lambda^*)\}$, $\lambda\in\Lambda$,
 where $\Lambda\subset\R^d$ is any countable set
 and $\{e_\lambda^*\}$ are arbitrary elements of the dual 
 space $(L^p(\Omega))^*$;
\item \label{nsf:ii}
The same conclusion holds for $ p > 2$ provided that $\Omega$ has nonempty interior.
\end{enumerate-num}
\end{thm}

Let us emphasize that we do not impose any a priori assumptions 
on the countable set of frequencies $\Lambda\subset\R^d$ 
or on the coefficient functionals
$e_\lambda^*\in (L^p(\Omega))^*$. 

We leave open the case where $p > 2$ and the set $\Omega$ has no interior points,
e.g.\  $\Omega\subset\R$ is a Cantor-type set of positive measure.

\subsection{}
In view of Theorems \ref{thm:bounded_unc_frames}
and \ref{thm:noframes_exp}, the situation is as follows:
the space $L^p(\Omega)$, $p \ne 2$, does not admit
any unconditional Schauder frame of exponentials
(at least, if $\Omega$ has nonempty interior),
but for $p>1$  it does admit an unconditional Schauder frame 
$\{(g_n, g_n^*)\}$
consisting of seminormalized and uniformly bounded 
(or even stronger, unimodular) functions $\{g_n\}$.

In this connection it is natural to ask whether in
\thmref{thm:bounded_unc_frames} one may
choose  the coefficient functionals  $\{g_n^*\}$ 
so that they are also seminormalized
(as elements of the dual space $(L^p(\Omega))^*$).
To motivate the question, observe for instance that  
if $\{u_n\}$  is a seminormalized unconditional basis in a Banach space $X$,
 then the unique biorthogonal system $\{u_n^*\}$ is seminormalized 
 (see \cite[Section 1.6, Theorem 3]{You01}).

Our next result provides an answer to the above question,
which turns out to depend on $p$.

\begin{thm}\label{thm:seminormalized}
If $1 < p \le 2$ then the coefficient functionals
$\{g_n^*\}$ in \thmref{thm:bounded_unc_frames}
 can be chosen so that they are   seminormalized. To the contrary, 
 for $p > 2$ there is no unconditional Schauder frame 
 $\{(g_n, g_n^*)\}$ in $L^p(\Omega)$ such that 
 the functions $\{g_n\}$ are uniformly bounded
and  both  systems $\{g_n\}$ and $\{g^*_n\}$ are seminormalized.
\end{thm}

The rest of the paper is devoted to the proofs of Theorems
\ref{thm:bounded_unc_frames}, \ref{thm:noframes_exp}
and \ref{thm:seminormalized}.


\section{Preliminaries}

In this section we review some basic facts that will be used later in the proofs.

\subsection{}
We start with a basic property of unconditional Schauder frames, that can be proved
using the uniform boundedness principle in the same way 
as a similar statement for unconditional Schauder bases (see e.g.\ \cite[Proposition 3.1.3]{AK16}).

\begin{lemma}
\label{lem:unc_constant}
Let $\{(u_j, u_j^*)\}_{j=1}^\infty$ be an unconditional Schauder frame in 
a Banach space $X$.  Then there exists a constant $K$ such that for any 
$x\in X$ and for any
sequence of scalars $\{\theta_j\}$ satisfying  $|\theta_j|\le 1$, we have
\begin{equation}
\label{eqtheta}
\Big\| \sum_{j=1}^\infty \theta_j u_j^*(x) u_j \Big\| \leq K \|x\|
\end{equation}
and the series in \eqref{eqtheta} converges unconditionally.
\end{lemma}

In this case we say that  
$\{(u_j, u_j^*)\}_{j=1}^\infty$ 
is a \emph{$K$-unconditional} Schauder frame.

\subsection{}
The following property of unconditional Schauder frames will also be useful.

\begin{lemma}\label{lem:changeroles}
Let $\{(u_j, u_j^*)\}_{j=1}^\infty$ be an unconditional Schauder frame in a reflexive Banach space $X$. Then the system $\{(u_j^*, u_j)\}_{j=1}^\infty$ is an unconditional Schauder frame in $X^*$.
\end{lemma}

See \cite[Section 4.3]{LT25} where a short proof of this fact is outlined.

Moreover, by duality it is easy
 to show  that if 
   $\{(u_j, u_j^*)\}_{j=1}^\infty$ 
 is  $K$-unconditional, then also
 $\{(u_j^*, u_j)\}_{j=1}^\infty$ 
 is $K$-unconditional.

\subsection{}
We use $\{r_n(t)\}_{n=1}^{\infty}$ to denote the sequence of Rademacher functions 
defined on the interval $[0,1]$ by $r_n(t) = \sign \sin (2^n \pi t)$. The classical Khintchine inequality states that for any $1\le p < +\infty$ there exist positive constants $A_p$ and $B_p$ such that for any $N$ and any sequence of scalars $\{a_n\}$ we have
\begin{equation}
\label{eq:KH.1}
A_p \Big( \sum_{n=1}^N |a_n|^2 \Big)^{1/2}\le \Big(\int_0^1 \Big| \sum_{n=1}^N a_nr_n(t) \Big|^p dt\Big)^{1/p}\le B_p  \Big( \sum_{n=1}^N |a_n|^2 \Big)^{1/2}
\end{equation}
(see \cite[I.B.8]{Woj91}). Note that
since the system $\{r_n\}_{n=1}^\infty$ is orthonormal in $L^2[0,1]$, one may take $B_p = 1$ for $p\le 2$, and $A_p = 1$ for $p\ge 2$.

The following lemma is a straightforward consequence of Khintchine's inequality.

\begin{lem}\label{lem:unc_conv_Rademachers}
The series $ \sum_{n=1}^\infty a_n r_n$
converges in $L^p[0,1]$,
$1\le p < +\infty$,
 if and only if $\sum_{n=1}^\infty |a_n|^2 < +\infty$. Moreover, in this case
 the series converges unconditionally and 
\begin{equation}
\label{eq:RD.5.1}
A_p \Big( \sum_{n=1}^\infty |a_n|^2 \Big)^{1/2}\le\Big\|\sum_{n=1}^\infty a_n r_n\Big\|_p\le  B_p \Big( \sum_{n=1}^\infty |a_n|^2 \Big)^{1/2}.
\end{equation}
\end{lem}

\subsection{}
The next lemma formulates the fact that  
 the space $L^p(\Om)$, $1\le p\le 2$, has type $p$ and cotype $2$.

\begin{lemma}
\label{lem:cotype}
Let $1\le p \le 2$. There exists a positive constant $C_p$ such that for any finite 
number of functions $f_1, f_2, \ldots, f_N \in L^p(\Omega)$ we have 
\begin{equation}
\label{eq:cotype}
C_p^{-1}\Big( \sum_{j=1}^N \|f_j\|_p^2 \Big)^{1/2}\leq \int_0^1 \Big\|\sum_{j=1}^N r_j(t) f_j\Big\|_p  dt\le C_p\Big( \sum_{j=1}^N \|f_j\|_p^p \Big)^{1/p}.
\end{equation}
\end{lemma}

A proof of \lemref{lem:cotype}, as well as other information about the notions of type and cotype of a Banach space, can be found e.g.\ in \cite[Chapter III.A]{Woj91}. 

 Lastly we establish the following estimate as a consequence of \lemref{lem:cotype}.

\begin{lemma}
\label{lem:averaging_Schframe}
Assume that $\{(g_j, g_j^*)\}_{j=1}^\infty$ is a $K$-unconditional Schauder frame in the space $L^p(\Omega)$, $1 < p\le 2$. Then for every $f \in L^p(\Om)$ we have
\begin{equation}
\label{eq:A.V.1}
\Big(\sum_{j=1}^\infty |g_j^*(f)|^2 \|g_j\|^2_p\Big)^{1/2} \le K C_p \|f\|_p.
\end{equation}
\end{lemma}

The proof is straightforward: due to \eqref{eqtheta},
for any $t\in [0,1]$ and every $N$ we have
\begin{equation}
\label{eq:A.V.2}
\Big\| \sum_{j=1}^{N}  r_j(t) g_j^*(f) g_j \Big\|_p \le K\|f\|_p,
\end{equation}
hence taking $\int_{0}^{1} dt$ of \eqref{eq:A.V.2},
 using \eqref{eq:cotype} and letting $N \to +\infty$,
 we arrive at \eqref{eq:A.V.1}.


\section{Unconditional Schauder frames of unimodular functions}
\label{sec:U1.1}

In this section we prove \thmref{thm:bounded_unc_frames},
i.e.\ we construct  an unconditional Schauder frame
$\{(g_n, g_n^*)\}$  in $L^p(\Omega)$, $1<p<+\infty$, consisting of
unimodular functions. 

We will prove two versions of the result. In the first version
the functions $g_n$ are \emph{real-valued} and unimodular,
i.e.\ they are $\{-1, +1\}$-valued. In the second version,
the functions $g_n$ are \emph{smooth} (complex-valued)  and unimodular.

\subsection{}
We begin with the first version, where 
the functions $g_n$  are $\{-1, +1\}$-valued.
With no loss of generality, we may assume that $\mes(\Omega) = 1$. Moreover, since all (separable, non-atomic) probability spaces are isomorphic (see \cite[\textsection 41]{Hal50}), we can in fact assume that $\Omega = [0,1]\subset \R$. We therefore need to prove the following assertion.

\begin{thm}\label{thm:uncframes[0,1]}
For any $p>1$ there exist in the space $L^p[0,1]$ a sequence $\{g_j\}_{j=0}^\infty$ 
consisting of $\{-1, +1\}$-valued functions, and a sequence $\{g_j^*\}_{j=0}^\infty$ of elements of the dual space $(L^p[0,1])^*$, such that any $f\in L^p[0,1]$ admits a series expansion
\begin{equation}
f=\sum_{j=0}^\infty g_j^*(f) g_j
\end{equation}
which converges unconditionally in $L^p[0,1]$.
\end{thm}

\subsection{}
Our proof is partly inspired by the paper \cite{FOSZ14}.
Following this paper we introduce the notion of an 
unconditional approximate Schauder frame.

\begin{definition}
	Let $X$ be a separable Banach space with dual space $X^*$. 
	A sequence of elements $\{(u_n, u_n^*)\}_{n=1}^\infty \subset X\times X^*$ is called an \emph{approximate Schauder frame} for $X$ if for any $x\in X$ the series
	\begin{equation}\label{eq:approx_frame}
		Sx = \sum_{n=1}^\infty u_n^*(x) u_n
	\end{equation}
	converges in $X$ and defines a bounded and invertible linear operator $S: X \to X$. If additionally the series \eqref{eq:approx_frame} converges unconditionally for every $x\in X$, then the system 
	$\{(u_n, u_n^*)\}_{n=1}^\infty$ is called an \emph{unconditional approximate Schauder frame}.
\end{definition} 

It is straightforward to verify that if $\{(u_n, u_n^*)\}_{n=1}^\infty$ is an 
unconditional approximate Schauder frame for $X$, 
then $\{(u_n, (S^{-1})^* u_n^*)\}_{n=1}^\infty$
 forms an unconditional Schauder frame for $X$ 
 (see \cite[Lemma 3.1]{FOSZ14}). 
 Therefore, in order to prove 
 \thmref{thm:uncframes[0,1]},
 it suffices to construct an 
unconditional approximate Schauder frame $\{(g_n, g_n^*)\}$
in $L^p[0,1]$  such that the functions $\{g_n\}$ are
$\{-1, +1\}$-valued.

In order to demonstrate our approach 
we will start by proving a slightly weaker statement:
we will construct an
unconditional approximate Schauder frame $\{(g_j, g_j^*)\}$
in $L^p[0,1]$  such that the functions $\{g_j\}$ are real-valued  and satisfy the condition
\begin{equation}
\label{eq:U.1}
1 - \eps \le |g_j(x)| \le 1 + \eps \quad \text{a.e.},
\end{equation}
where $\eps$ is any fixed positive number  given in advance.
Then we will show how to modify the construction
in order to obtain $\{-1, +1\}$-valued functions $\{g_j\}$.

\subsection{}	
We consider   the system of Haar functions 
$\{h_k\}_{k=0}^\infty$ normalized in the $L^p$ norm (we use the convention   $h_0\equiv 1$)
and denote the biorthogonal system by $\{h_k^*\}_{k=0}^\infty$. 
If we identify the dual space $(L^p[0,1])^*$ with $L^q[0,1]$, $q = p/(p-1)$, 
 in the usual way, then
the elements $h_k^*$ are  the Haar functions normalized in $L^q$. 
Hence $\|h_k^*\| = 1$ for all $k$.

We now fix a sequence of positive numbers $\{\del_k\}_{k=1}^\infty$ satisfying
\begin{equation}
\label{eq:D.K.1}
\Big(\sum_{k=1}^{\infty} \del_k^2\Big)^{1/2} < \del,
\end{equation}
where $\del = \del(p) > 0$ is a small number to be specified later.
Next we choose a sequence of large positive integers $\{N_k\}_{k=1}^\infty$ and consider a partition of the set $\{1,2,3,\dots\}$ into consecutive ``blocks'' $\{J_k\}_{k=1}^\infty$ such that $\# J_k = N_k$. We also recall that $\{r_j\}_{j=1}^{\infty}$ denotes the sequence of Rademacher functions on the interval $[0,1]$. 

Define $g_0 = h_0$, $g_0^* = h_0^*$, and
\begin{equation}\label{eq:uncframe_constuction}
g_j = \del_k^{-1} N_k^{-1/2} h_k +  r_j, \quad g_j^* = \del_k N_k^{-1/2}h_k^*,
\quad  j\in J_k.
\end{equation}
We claim that the system $\{(g_j, g_j^*)\}_{j=0}^\infty$ is an 
unconditional approximate Schauder frame  for the space $L^p[0,1]$ provided that $\delta = \del(p) $ is sufficiently small. Moreover, since we have $|r_j(x)| = 1$ a.e., we can clearly choose the numbers $N_k$ large enough so that the condition \eqref{eq:U.1} is satisfied, where $\eps>0$ is any fixed number given in advance.

We now need to verify that  for every $f\in L^p[0,1]$ the series
\begin{equation}\label{eq:to_check_approxframe}
Sf = \sum_{j=0}^\infty g_j^*(f) g_j
\end{equation}
converges unconditionally in $L^p[0,1]$ and defines a bounded and invertible operator $S$ on the space $L^p[0,1]$. 

If we put $J_0 = \{0\}$ and $N_0 = 1$, then by \eqref{eq:uncframe_constuction},
\eqref{eq:to_check_approxframe}  we have, at least formally,
\begin{equation}
\label{eq:S.1}
Sf = \sum_{k=0}^\infty \sum_{j\in J_k} N_k^{-1} h_k^*(f) h_k + \sum_{k=1}^\infty \sum_{j\in J_k} \del_k N_k^{-1/2} h_k^*(f) r_j,
\end{equation}
so it would suffice to treat each of the two series in \eqref{eq:S.1} separately. 
Since $\{h_k\}$ forms an unconditional basis in the space $L^p[0,1]$, and since
 $\# J_k = N_k$, it is easy to verify
(for instance, using  \cite[Propositions 2.4.9 and 3.1.3]{AK16})
  that the first series converges unconditionally in $L^p[0,1]$ to the function $f$.

Next, we claim that also the ``error term'' series 
\begin{equation}
\label{eq:E1.1}
Ef =  \sum_{k=1}^\infty \sum_{j\in J_k} \del_k N_k^{-1/2} h_k^*(f) r_j
\end{equation}
converges unconditionally in $L^p[0,1]$, and moreover,
 if $\del = \del(p)$  is sufficiently small then 
the norm of the operator $E$ defined by \eqref{eq:E1.1} is smaller than $1$. 
This will suffice in order to conclude the proof, since then the operator $S=I + E$ is invertible.

Indeed, by \lemref{lem:unc_conv_Rademachers} 	the
series \eqref{eq:E1.1} converges unconditionally provided that
\begin{equation}
\label{eq:E1.3}
\Big(\sum_{k=1}^\infty \del_k^2 \, |h_k^*(f)|^2\Big)^{1/2}
\end{equation}
is finite. 
Recalling that $\|h_k^*\|= 1$ and using \eqref{eq:D.K.1} it follows that 
\eqref{eq:E1.3} does not exceed
\begin{equation}
\Big( \sum_{k=1}^\infty \del_k^2  \Big)^{1/2} \|f\|_p \le \del \|f\|_p.
\end{equation}
Hence the
series \eqref{eq:E1.1} converges unconditionally,
and moreover $\|E f\|_p \le B_p \del \|f\|_p$.
So by choosing $\delta$ such that $B_p \del < 1$
we can ensure that
the norm of the operator $E$ is indeed smaller than $1$. 
Thus the system $\{(g_j, g_j^*)\}_{j=0}^{\infty}$   forms an 
unconditional approximate Schauder frame in 
the space $L^p[0,1]$ as we had to show.

\subsection{}
Next, we show how to modify the above construction so as to obtain real-valued
functions $\{g_j\}$ such that $|g_j(x)| =1$ a.e., that is,
$g_j$ are $\{-1,+1\}$-valued functions.
To  this end, we choose the sequence $\{\del_k\}_{k=1}^\infty$ to satisfy 
not the condition \eqref{eq:D.K.1} but rather the stronger condition
\begin{equation}
\label{eq:3.10}
\sum_{k=1}^{\infty} \del_k < \del.
\end{equation} 

Now, recalling that for each $k$ the function $|h_k(x)|$ is constant on its support,
we may choose $N_k$ large enough and adjust the parameter $\del_k$ 
so that $\del_k^{-1} N_k^{-1/2} |h_k(x)| = 1$ for
$x \in \supp(h_k)$. Then, in the definition
\eqref{eq:uncframe_constuction} we modify the function $g_j$ 
and choose it as
\begin{equation}
\label{eq:unimod_def}
g_j(x) = \del_k^{-1} N_k^{-1/2} h_k(x) + (1- \chi_k(x)) r_j(x), 
\quad  j\in J_k,
\end{equation}
where $\chi_k$ denotes the indicator function of $\supp(h_k)$.
Then   $g_j$  is a real-valued function with $|g_j(x)| = 1$ a.e. 
We keep the same coefficient functionals $g_j^*$ as before, namely
\begin{equation}
\label{eq:unimod_coeff}
g_j^* = \del_k N_k^{-1/2}h_k^*,	\quad  j\in J_k.
\end{equation}

We now observe that for each  $k$ the estimate 
\begin{equation}
\label{eq:E5.1}
\Big\| \sum_{j\in J_k} \del_k N_k^{-1/2} h_k^*(f) (1-\chi_k(x)) c_j r_j(x) \Big\|_p
\le B_p \del_k \|f\|_p
\end{equation}
holds for any sequence  of scalars $\{c_j\}$ with 
$|c_j| \le 1$. Indeed, this follows from
\lemref{lem:unc_conv_Rademachers} using the
fact that $\|h_k^*\|=1$, and that
multiplication by the function $1-\chi_k$ does not increase the $L^p$ norm.

Due to the condition \eqref{eq:3.10}, the estimate \eqref{eq:E5.1}
implies (for instance, using \cite[Proposition 2.4.9]{AK16})
that the modified ``error term'' series
\begin{equation}
	\label{eq:E5.3}
	(Ef)(x) =  \sum_{k=1}^\infty \sum_{j\in J_k} \del_k N_k^{-1/2} h_k^*(f) (1-\chi_k(x)) r_j(x)
\end{equation}
converges unconditionally in $L^p[0,1]$, and satisfies
$\|E f\|_p \le B_p \del \|f\|_p$. As before,
this suffices to ensure that the system $\{(g_j, g_j^*)\}_{j=0}^{\infty}$
forms an unconditional approximate Schauder frame in $L^p[0,1]$.
Thus the full statement of 
 \thmref{thm:uncframes[0,1]} is now proved. \qed

\subsection{}
Finally we also prove the second version of
\thmref{thm:bounded_unc_frames}.
 Let $\Omega\subset\R^d$ be a set of positive and finite measure.
For any $p>1$ we will now construct an unconditional Schauder frame 
$\{(\pphi_n, \pphi_n^*)\}$ in the space $L^p(\Omega)$ 
such that the functions $\pphi_n$ are
uni\-modular (complex-valued) \emph{and smooth},
i.e.\ each $\pphi_n$ is obtained as the restriction
to $\Om$ of a unimodular and smooth function in the whole $\R^d$.

Indeed, let $\{(g_n, g_n^*)\}$ be the
unconditional Schauder frame constructed above
for the space $L^p(\Om)$. Since $|g_n(x)|=1$ a.e.,
we can write $g_n(x) = \exp 2 \pi i u_n(x)$
 where $u_n$ is a real-valued function on $\Om$.
Given $\eps_n>0$ we can find a real-valued and smooth function
 $v_n$ on $\R^d$ which approximates $u_n$ in measure on $\Om$,
 that is,
\begin{equation}
\label{eq:SA.1.1}
\mes \{x \in \Om : |v_n(x)-u_n(x)| > c   \eps_n \} <  (c \eps_n)^p,
\end{equation}
where $c = c(p,\Om)>0$ is a small constant. 
Let $\pphi_n(x) = \exp 2 \pi i v_n(x)$,
then $\pphi_n$ is a smooth function on $\R^d$
with $|\pphi_n(x)|=1$, and \eqref{eq:SA.1.1} implies that 
$\|\pphi_n - g_n\|_{L^p(\Om)} < \eps_n$.
If we choose $\eps_n$ so that
$\sum_n \eps_n \|g^*_n\| < 1$, then $\{(\pphi_n, g_n^*)\}$ 
forms an unconditional approximate Schauder frame in 
$L^p(\Om)$. Finally, by an application of
\cite[Lemma 3.1]{FOSZ14} we obtain
functionals $\{\pphi^*_n\}$ such that
$\{(\pphi_n, \pphi_n^*)\}$ is
an unconditional Schauder frame  in $L^p(\Omega)$.

\subsection{Remark}
\label{sec:R3.1}

One may ask whether there exists an unconditional Schauder frame 
in the space $L^p(\Omega)$ consisting of \emph{nonnegative functions}. 
The answer is negative \cite[Theorem 5.5]{PS16}, see also \cite[Section 3]{JS15}.
In this connection, we also mention that there do exist
(not unconditional) Schauder frames 
in $L^p(\Om)$ formed by nonnegative functions 
  \cite[Corollary 5.4]{PS16}.
Moreover, in the spaces $L^1(\Om)$ and $L^2(\Om)$
there even exists 
a Schauder basis consisting of nonnegative functions,
 see  \cite{JS15}, \cite{FPT21}.
The existence of such a basis in $L^p(\Om)$, $p\neq 1, 2$, remains open.


\section{The existence of seminormalized coefficient functionals}

Now we turn to the proof of \thmref{thm:seminormalized}.
We will show
that if $1 < p \le 2$ then the coefficient functionals  in 
\thmref{thm:bounded_unc_frames} 
 can be chosen so that they are seminormalized.
 To the contrary, for $p>2$ we will show that if
 $\{(g_n, g_n^*)\}$ 
 is any unconditional Schauder frame in $L^p(\Omega)$, and if
 the functions $\{g_n\}$ are seminormalized in $L^p(\Om)$ and uniformly bounded,
 then the coefficient functionals  $\{g^*_n\}$ cannot be   seminormalized.

As in \secref{sec:U1.1}, we may assume that $\Omega = [0,1]\subset \R$.

\subsection{The case $1<p\le 2$}
We start with the first part of the theorem.
 Recall that in the construction above the system $\{(g_j, g_j^*)\}$ 
 defined by \eqref{eq:unimod_def}, \eqref{eq:unimod_coeff}
forms an unconditional approximate Schauder frame in $L^p[0,1]$,
 which in turn implies that $\{(g_j, (S^{-1})^* g_j^*)\}$ 
 is an unconditional Schauder frame, where $S$ is the (bounded, invertible)
 operator defined by \eqref{eq:to_check_approxframe}.
However, note that since 
the numbers  $\delta_k$ satisfy \eqref{eq:3.10}, and moreover,
there are arbitrarily large $N_k$'s, neither the system 
$\{  g_j^* \}$ nor $\{ (S^{-1})^* g_j^* \}$ is seminormalized.

We will nevertheless show that if $1 < p \le 2$ then
these  coefficient functionals can be 
``corrected'' so that they become seminormalized. 
This will be done by an application of 
\cite[Corollary 2.5]{BC20} which can be stated as follows:

\begin{lem}
\label{lem:BC}
Let $\{(g_j, g_j^*)\}_{j=0}^\infty$ be an
 approximate Schauder frame for a Banach space $X$ such that $\{g_j\}_{j=0}^\infty$ are seminormalized. Suppose that $X$ contains a complemented copy of $\ell^r$ for some $1\le r < +\infty$, and that there exists a positive constant $M$ such that
\begin{equation}
\label{eq:BC.1.2}
	\Big\| \sum_{j=0}^N a_j g_j \Big\| \le M \Big( \sum_{j=0}^N |a_j|^r \Big)^{1/r}
\end{equation}
for every $N$ and every sequence of scalars $\{a_j\}$. Then there exists a seminormalized sequence $\{G_j^*\}_{j=0}^\infty \subset X^*$ such that $\{(g_j, G_j^*)\}_{j=0}^\infty$ is a Schauder frame for $X$. Moreover, if $\{(g_j, g_j^*)\}_{j=0}^\infty$ is an unconditional approximate Schauder frame, then $\{G_j^*\}_{j=0}^\infty$ can be chosen so that $\{(g_j, G_j^*)\}_{j=0}^\infty$ forms an unconditional Schauder frame.
\end{lem}

We wish to apply this lemma in the space $X=L^p[0,1]$, $1 < p \le 2$,  with $r=p$ and  the
unconditional approximate Schauder frame
$\{(g_j, g_j^*)\}_{j=0}^\infty$ given by \eqref{eq:unimod_def},
 \eqref{eq:unimod_coeff}.
It is well-known and easy to prove that $L^p[0,1]$ contains a complemented 
subspace isometric to $\ell^p$ (see e.g.\ \cite[Lemma 1]{KP62}). 
However, an obstacle appears when we try to establish
the estimate \eqref{eq:BC.1.2}, due to the large factors $\del_k^{-1}$ which appear
in the definition of $g_j$. 

To resolve this we shall now modify the sequence 
$\{\del_k\}_{k=1}^\infty$ and choose it to be an ``almost constant''
sequence. First we fix a small number $\delta=\delta(p) > 0$
that will be specified later. Then 
 we choose the parameters $\del_k$ and  $N_k$ so that
 \begin{equation}
\label{eq:parameters_seminormalized_1}
\delta_k^{-1} N_k^{-1/2} |h_k(x)| = 1, \quad x \in \supp(h_k),
\end{equation} 
and additionally, 
 \begin{equation}
\label{eq:parameters_seminormalized_2}
\tfrac1{2} \delta \le\delta_k\le\delta.
\end{equation} 
To see that such a choice of the parameters indeed exists, it suffices
to let $N_k$  be the smallest positive integer satisfying
$\delta^{-1} N_k^{-1/2} \|h_k\|_\infty \le 1$, 
then  define 
$\del_k$ by the condition \eqref{eq:parameters_seminormalized_1},
and finally check that \eqref{eq:parameters_seminormalized_2} holds.

However, note that in this case the condition \eqref{eq:3.10} is clearly no longer satisfied, 
which in turn means that our estimate of the ``error term'' series
\eqref{eq:E5.3} is no longer valid. Nevertheless we can take advantage of
the assumption that $1<p\le 2$ and still obtain the required estimate by
a more careful analysis.

\subsubsection{The estimation of the error term}

 Let $K_p$ denote the unconditionality constant of the Haar basis in $L^p[0,1]$,
 which by definition is the least constant $K_p$ such that
\begin{equation}
\label{eq:H.4.1}
\Big\| \sum_{k=0}^{N} b_k h_k \Big\|_p \le K_p
\Big\| \sum_{k=0}^{N} a_k h_k  \Big\|_p
\end{equation}
holds for any $N$ and any two sequences of scalars $\{a_k\}, \{b_k\}$
satisfying $|b_k| \le |a_k|$.
(In fact, it is known that $K_p = (p-1)^{-1}$ for $1<p\le 2$, see \cite[Theorem II.D.13]{Woj91}.)

We recall that our system $\{(g_j, g_j^*)\}$ is still given by \eqref{eq:unimod_def},
 \eqref{eq:unimod_coeff},
 but this time the parameters $\delta_k$ and $N_k$ satisfy the conditions 
 \eqref{eq:parameters_seminormalized_1} and  \eqref{eq:parameters_seminormalized_2}.
  Therefore, the error term which we need to 
   estimate is still given by \eqref{eq:E5.3}, which we decompose as
\begin{equation}
\label{eq:to_est_seminorm}
(Ef)(x) = \sum_{k=1}^\infty \sum_{j\in J_k} 
\del_k N_k^{-1/2} h_k^*(f)  r_j(x) 
- \sum_{k=1}^\infty \sum_{j\in J_k} 
\del_k N_k^{-1/2} h_k^*(f) \chi_k(x) r_j(x).
\end{equation}

The first series in \eqref{eq:to_est_seminorm}
converges unconditionally due to \lemref{lem:unc_conv_Rademachers},
which also provides the estimate
\begin{equation}
	\Big\|\sum_{k=1}^\infty \sum_{j\in J_k} 
	\del_k N_k^{-1/2} h_k^*(f)  r_j\Big\|_p
	\le  B_p \delta \Big( \sum_{k=1}^\infty |h_k^*(f)|^2 \Big)^{1/2},
\end{equation}
and by an application of \lemref{lem:averaging_Schframe} to the normalized Haar basis we get that this quantity does not exceed $\delta B_pK_pC_p\|f\|_p$.

It remains to estimate the second series in \eqref{eq:to_est_seminorm} and to establish that it converges unconditionally. In order to do it, we recall that according to the standard enumeration of the Haar functions (see  \cite[II.B.9]{Woj91}), for any fixed integer $l\ge 0$
the supports of the Haar functions $\{ h_k \}$, $2^l\le k<2^{l+1}$, consist of
 disjoint dyadic intervals of length $2^{-l}$. Recalling that $\chi_k$ is the indicator 
 function of $\supp(h_k)$, it follows that
\begin{align}
& \Big\|\sum_{2^l\le k < 2^{l+1}} \sum_{j\in J_k} \del_k N_k^{-1/2} h_k^*(f)  c_j  \chi_k(x) r_j(x)\Big\|_p^p \\
& \qquad
	 = \sum_{2^l\le k < 2^{l+1}} \Big\|  \chi_k(x)   \sum_{j\in J_k} \del_k N_k^{-1/2} h_k^*(f)  c_j r_j(x)  \Big\|_p^p \label{eq:4.8}
\end{align}
for any sequence of scalars $\{c_j\}$. Now observe that each
$r_j$ is a $2^{-j+1}$-periodic  function, hence 
if $k \ge 2^l$ then all the functions $\{r_j\}$, $j \in J_k$, are $2^{-l}$-periodic.
Since $\chi_k$ is the indicator function of a dyadic interval of length $2^{-l}$,
the quantity \eqref{eq:4.8} is thus equal to
\begin{equation}
\label{eq:D.2.1}
\sum_{2^l\le k < 2^{l+1}}  2^{-l} \, \Big\|\sum_{j\in J_k} \del_k N_k^{-1/2} h_k^*(f)  c_j r_j(x)  \Big\|_p^p.
\end{equation}
Now  suppose that $|c_j|\le 1$, then by \lemref{lem:unc_conv_Rademachers} 
the quantity \eqref{eq:D.2.1} does not exceed
\begin{equation}
B_p^p \delta^p  2^{-l}  \sum_{2^l\le k < 2^{l+1}} |h_k^*(f)|^p.
\end{equation}
We have therefore shown that for each   $l \ge 0$  and
any scalars $\{c_j\}$ with $|c_j| \le 1$, we have
\begin{equation}
\Big\|\sum_{2^l\le k < 2^{l+1}}  \sum_{j\in J_k} \del_k N_k^{-1/2} h_k^*(f)  
c_j \chi_k(x) r_j(x)\Big\|_p
\le \delta B_p \Big( 2^{-l}\sum_{2^l\le k < 2^{l+1}} |h_k^*(f)|^p \Big)^{1/p}.
\end{equation}
Since $1 < p\le 2$, we can estimate the latter quantity by
\begin{equation}
\label{eq:S.4.1}
\delta B_p \Big( 2^{-l}\sum_{2^l\le k < 2^{l+1}}
 |h_k^*(f)|^2 \Big)^{1/2}\le \delta B_p K_pC_p 2^{-l/2} \|f\|_p,
\end{equation}
where in the last inequality we once again applied \lemref{lem:averaging_Schframe} to the normalized Haar basis. It remains to use the triangle inequality in $L^p$ and conclude that
the second series in \eqref{eq:to_est_seminorm} converges unconditionally
and satisfies the estimate 
\begin{equation}
\label{eq:S.4.2}
\Big\|\sum_{k=1}^\infty \sum_{j\in J_k}
 \del_k N_k^{-1/2} h_k^*(f) 
 \chi_k(x) r_j(x)\Big\|_p\le (2+\sqrt{2})
  \delta B_p K_pC_p \|f\|_p,
\end{equation}
where the constant $2 + \sqrt{2}$ comes from summing the terms
$2^{-l/2}$ in \eqref{eq:S.4.1}.

We conclude that if $\delta = \del(p) > 0$ is small enough, then 
the norm of the operator $E$ given by \eqref{eq:to_est_seminorm}
 is indeed smaller than $1$, hence the
  system $\{(g_j, g_j^*)\}_{j=0}^{\infty}$
  defined by \eqref{eq:unimod_def},  \eqref{eq:unimod_coeff}
  (as well as $g_0 = h_0$, $g_0^* = h_0^*$)   
with parameters $\delta_k$ and $N_k$ satisfying
 \eqref{eq:parameters_seminormalized_1},
  \eqref{eq:parameters_seminormalized_2}
 forms an unconditional approximate Schauder frame for the space $L^p[0,1]$.

\subsubsection{An application of \lemref{lem:BC}}
It remains now to check that \lemref{lem:BC} 
(with $r=p$) is applicable to the system $\{(g_j, g^*_j)\}$. 
That is, we need to show that  for arbitrary scalars 
$\{a_j\}_{j=1}^\infty$ such that only finitely 
many of them are nonzero, the inequality
\begin{equation}
\label{eq:H.J.12}
\Big\| \sum_{k=0}^\infty \sum_{j\in J_k} a_j g_j \Big\|_p \le M \Big( \sum_{k=0}^\infty \sum_{j\in J_k} |a_j|^p \Big)^{1/p}
\end{equation}
holds for some constant $M$. By the definition \eqref{eq:unimod_def}
 of $g_j$ we have
\begin{align}
& \Big\| \sum_{k=1}^\infty \sum_{j\in J_k} a_j g_j \Big\|_p = \Big\| \sum_{k=1}^\infty \sum_{j\in J_k} a_j ( \del_k^{-1} N_k^{-1/2} h_k +  (1-\chi_k)r_j) \Big\|_p \label{eq:BC.2.1} \\
& \qquad
\le \Big\| \sum_{k=1}^\infty \Big(\delta_k^{-1} N_k^{-1/2} \sum_{j\in J_k} a_j \Big) h_k \Big\|_p + \Big\| \sum_{j=1}^\infty  a_j r_j \Big\|_p+ \Big\| \sum_{k=1}^\infty \sum_{j\in J_k} a_j \chi_k r_j \Big\|_p. \label{eq:BC.2.2}
\end{align}
The second summand in \eqref{eq:BC.2.2} does not exceed
\begin{equation}
 \Big\| \sum_{j=1}^\infty a_j r_j \Big\|_2 =  \Big( \sum_{j=1}^\infty |a_j|^2 \Big)^{1/2}\le  \Big( \sum_{j=1}^\infty |a_j|^p \Big)^{1/p}.
\end{equation}

The third summand can be estimated in a similar way: it suffices to
observe that  the functions $\{\chi_k r_j\}$, $k \ge 1$, $j \in J_k$,
  are pairwise orthogonal in $L^2[0,1]$. 
 Indeed, $\supp(\chi_k)$ is a dyadic interval of length $2^{-l}$ for 
 $2^l\le k < 2^{l+1}$, and since $j\ge k$, $r_j$ is non-constant on 
 $\supp(\chi_k)$. It remains to observe that if at least one of two distinct 
 Rademacher functions $r_{j}$ and $r_{j'}$ is non-constant on a 
 dyadic interval $I$, then these functions are orthogonal in $L^2(I)$.
 Hence the third summand in \eqref{eq:BC.2.2} does not exceed
\begin{equation}
\Big\| \sum_{k=1}^\infty \sum_{j\in J_k} a_j \chi_k r_j \Big\|_2
\le \Big( \sum_{j=1}^\infty |a_j|^2 \Big)^{1/2}
\le  \Big( \sum_{j=1}^\infty |a_j|^p \Big)^{1/p}.
\end{equation}

In order to estimate the first summand in \eqref{eq:BC.2.2} 
 we will use the fact that for any $N$ and for any sequence of scalars
  $\{c_k\}$ we have
\begin{equation}
\label{eq:H.C.5}
\Big\| \sum_{k=0}^N c_k h_k \Big\|_p \le K_p C_p \Big(\sum_{k=0}^N |c_k|^p\Big)^{1/p}.
\end{equation}
Indeed, \eqref{eq:H.C.5} can be proved similarly to 
 \lemref{lem:averaging_Schframe}: by \eqref{eq:H.4.1}
for any $t\in [0,1]$  we have
\begin{equation}
\label{eq:H.C.6}
\Big\| \sum_{k=0}^N c_k h_k \Big\|_p \le
K_p \Big\| \sum_{k=0}^N r_k(t) c_k h_k \Big\|_p,
\end{equation}
so taking $\int_{0}^{1} dt$ of \eqref{eq:H.C.6} and
applying \lemref{lem:cotype} we arrive at \eqref{eq:H.C.5}.

As a consequence of \eqref{eq:H.C.5}, 
the first summand in \eqref{eq:BC.2.2}  does not exceed
\begin{equation}\label{eq:3.15}
K_pC_p \Big(\sum_{k=1}^\infty \delta_k^{-p} N_k^{-p/2}
  \Big| \sum_{j\in J_k} a_j \Big|^p \Big)^{1/p}.
\end{equation}
By H\"{o}lder's inequality, we have
\begin{equation}
\label{eq:H2.1}
\Big| \sum_{j\in J_k} a_j \Big|^p\le N_k^{p-1} \sum_{j\in J_k}|a_j|^p.
\end{equation}
Hence \eqref{eq:parameters_seminormalized_2} and \eqref{eq:H2.1} 
imply that \eqref{eq:3.15} is not greater than
\begin{equation}
2\del^{-1} K_pC_p\Big( \sum_{k=1}^\infty N_k^{(p-2)/2}
\sum_{j\in J_k}|a_j|^p \Big)^{1/p}
\le 2\del^{-1} K_pC_p \Big( \sum_{k=1}^\infty
\sum_{j\in J_k} |a_j|^p \Big)^{1/p},
\end{equation}
where  the last inequality is true since we have
$(p-2)/2 \le 0$ and $N_k\ge 1$. The required estimate \eqref{eq:H.J.12}
thus follows.
To conclude, all conditions of \lemref{lem:BC} have been checked and  
the proof of the first part of \thmref{thm:seminormalized} is therefore 
complete. \qed

\subsection{The case $p > 2$}	
Now we prove the second part of \thmref{thm:seminormalized}. 
It suffices to prove the following assertion.

\begin{thm}
\label{thm:S.1}
Suppose that $\{(g_j, g_j^*)\}_{j=1}^\infty$ is
an unconditional Schauder frame for the space 
$L^p[0,1]$, $p>2$, such that
the functions  $\{g_j\}$ are
seminormalized in $L^p[0,1]$ and
they are uniformly bounded. Then 
\begin{equation}
\label{eq:M.P.2}
\inf_j \|g_j^*\| =0,
\end{equation}
and as a consequence, the system
 $\{g^*_j\}$ cannot be seminormalized.
\end{thm}

Our proof is partly inspired by the paper \cite{KP62}.

\begin{proof}
Let $M$ be a constant  such that for every $j$ we have
 $|g_j(x)|\le M$ a.e. With no loss of generality we 
may assume that $\|g_j\|_p = 1$.  In this case we have
\begin{equation}
1= \|g_j\|_p^p 
\le \|g_j\|_\infty^{p-2} \|g_j\|_2^{2}\le M^{p-2}\|g_j\|_2^{2},
\end{equation}
 hence $\|g_j\|_2\ge  M^{(2-p)/2}$. 
The functions $\{g_j\}$ are thus seminormalized  in  $L^2[0,1]$
 (in fact, they are seminormalized in all $L^r$ norms, $1\le r\le +\infty$, 
 but we do not use this).

By \lemref{lem:unc_constant} there is a constant $K$ such that 
\begin{equation}
\Big\|\sum_{j=1}^{N} r_j(t)g_j^*(f)g_j\Big\|_2^2
\le 	\Big\|\sum_{j=1}^{N} r_j(t)g_j^*(f)g_j\Big\|_p^2
\le K^2\|f\|_p^2
\end{equation}
holds for every $f\in L^p[0,1]$, every $N$ and every $t\in [0,1]$.
Since the Rademacher functions 
$\{r_j\}$ are orthonormal in the space $L^2[0,1]$, then
by taking $\int_0^1 dt$,  applying the previously obtained 
estimate for $\|g_j\|_2$, and letting $N \to +\infty$, we obtain
\begin{equation}
\label{eq:l2estem}
\Big(\sum_{j=1}^\infty |g_j^*(f)|^2\Big)^{1/2} \le M^{(p-2)/2}K\|f\|_p.
\end{equation}

Suppose to the contrary that \eqref{eq:M.P.2} does not hold.
Then there is $\del > 0$ such that $\|g^*_j\| \ge\delta$ for all $j$.
We now invoke \lemref{lem:changeroles}, which yields that
 the system $\{(g_j^*, g_j)\}_{j=1}^\infty$ forms an unconditional Schauder 
 frame for the dual space $(L^p[0,1])^* = L^{q}[0,1]$, where $q = p/(p-1)$.
Moreover, $\{(g_j^*, g_j)\}_{j=1}^\infty$ is
 a $K$-unconditional Schauder frame. In turn, by applying
 \lemref{lem:averaging_Schframe} (note that $1 < q < 2$)
  together with the fact that $\|g^*_j\|_{q}\ge\delta$, 
  we obtain that for any element $f^*\in (L^p[0,1])^*$ we have
\begin{equation}\label{eq:l2estem_reverse}
\Big(\sum_{j=1}^\infty |f^*(g_j)|^2\Big)^{1/2}\le \delta^{-1} K C_q   \|f^*\|_{q}.
\end{equation}
We can now use a standard duality argument. Since for any
two elements $f\in L^p[0,1]$ and $f^*\in (L^p[0,1])^*$ 
we have the unconditionally convergent expansions
\begin{equation}
f = \sum_{j=1}^\infty g_j^*(f) g_j,\quad f^* = \sum_{j=1}^\infty f^*(g_j) g_j^*,
\end{equation}
we have the estimate
\begin{equation}
\label{eq:A.6.1}
|f^*(f)| = \Big|\sum_{j=1}^\infty f^*(g_j) g_j^*(f)\Big|
\le 	\Big(\sum_{j=1}^\infty |f^*(g_j)|^2\Big)^{1/2} 
\Big(\sum_{j=1}^\infty |g_j^*(f)|^2\Big)^{1/2}.
\end{equation}
If we choose $f^*$ so that $\|f^*\|_{q} = 1$, $f^*(f) = \|f\|_p$,
and use \eqref{eq:l2estem_reverse}, \eqref{eq:A.6.1} then we obtain
\begin{equation}
\label{eq:A.6.2}
\|f\|_p \le \delta^{-1} K C_q  
\Big(\sum_{j=1}^\infty |g_j^*(f)|^2\Big)^{1/2}.
\end{equation}

The two inequalities \eqref{eq:l2estem}, \eqref{eq:A.6.2}
hold for any $f\in L^p[0,1]$. This means that the mapping
$Tf = \{g_j^*(f)\}_{j=1}^\infty$ is a bounded linear operator 
from $L^p[0,1]$ to $\ell^2$, and it is an isomorphism 
onto its (closed) image. Hence we get that the space $L^p[0,1]$
 is isomorphic to a closed subspace of $\ell^2$, namely, to a Hilbert space, which is clearly not true
(e.g.\ due to \cite[Theorem 6.2.14]{AK16}).
This contradiction thus completes the proof. 
\end{proof}

\subsection{Remark}
The assumption that 
the functions $\{g_j\}$ are uniformly bounded
was used in the proof above only in order
to conclude that $\inf_j \|g_j\|_2 > 0$.
Hence \thmref{thm:S.1} remains true
if the uniform boundedness assumption is replaced with
e.g.\ the weaker assumption that
$\sup_j \|g_j\|_r < + \infty$ for some $r > p$.


\section{Unconditional Schauder frames of exponentials: the case $1 < p < 2$}

In the present section we prove part \ref{nsf:i} of \thmref{thm:noframes_exp}.

\subsection{}
Suppose that there exists an unconditional Schauder frame
 $\{(e_\lambda, e_{\lambda}^*)\}$, $\lambda\in\Lambda$,
  in the space $L^p(\Omega)$, $1 < p < 2$, for some set $\Omega\subset\R^d$ of finite measure. Then the same system also forms an unconditional Schauder frame in $L^p(\Omega')$ for any $\Omega'\subset\Omega$. Hence we may assume with no loss of generality that the set $\Omega$ is bounded. Moreover, after rescaling we can also assume that $\mes(\Omega) = 1$. So we need to prove the following:

\begin{thm}
Let $\Omega \sbt \R^d$ be a bounded measurable set,  $\mes(\Omega) = 1$. Then 
in the space $L^p(\Omega)$, $1 < p < 2$, 
there is no unconditional Schauder frame 
of  the form $\{(e_\lambda, e_{\lambda}^*)\}$,  $\lambda\in\Lambda$,
 where $\Lambda\subset\R^d$ is a countable set and $\{e_\lambda^*\}$ are arbitrary elements of $(L^p(\Omega))^*$. 
\end{thm}

Note that we exclude the case $p=1$ since, as we mentioned,
it is known that there are no unconditional Schauder frames (of any form)
in the space $L^1(\Omega)$.

The proof below would be simpler if the set of frequencies
 $\Lambda$ was assumed to be  uniformly discrete.
In order to establish the result 
in the general case, i.e.\  for an arbitrary countable set $\Lambda
 \sbt \R^d$, we use an approach from \cite{LT25}.

\subsection{}
Suppose to the contrary  that $\{(e_\lambda, e_{\lambda}^*)\}$,
$\lambda\in\Lambda$, forms
 an unconditional Schauder frame in $L^p(\Omega)$.
 Our goal is to show that this assumption leads to a contradiction.
 
Since $\Omega$ is a bounded set, there exists $ \delta = \del(\Om) > 0$ with the following property: for every $\lambda, \lambda'\in \R^d$ such that $|\lambda'-\lambda| \le \delta$ we have $\|e_{\lambda'} - e_{\lambda}\|_{L^{p}(\Omega)} < 1/2$. We fix this number $\delta$ and 
consider a partition of $\R^d$ into congruent (half-open) cubes 
$\{Q_k\}_{k=1}^\infty$ with $\diam (Q_k) = \delta$. 
It induces a partition of $\Lam$ into ``blocks''  $\Lam_k =  \Lambda\cap Q_k$. 
We note that each block $\Lam_k$ may contain arbitrarily many
elements, or  may  even be infinite.

Now, every $f\in L^p(\Omega)$ has 
an unconditionally convergent series expansion
\begin{equation}
\label{eq:frame_expansion}
	f=\sum_{k=1}^\infty \sum_{\lambda\in\Lambda_k} e_\lambda^*(f) e_\lambda.
\end{equation}
By \lemref{lem:unc_constant} there exists a constant $K$ (not depending on $f$) such that 
\begin{equation}
	\Big\| \sum_{k=1}^{N} r_k(t) \sum_{\lambda\in{\Lambda_k}} |e_\lambda^*(f)| e_\lambda \Big\|_p \le K \|f\|_p
\end{equation}
holds for every $N$ and every $t\in [0,1]$. If we now take $\int_0^1 dt$,
 apply \lemref{lem:cotype},  and then let $N \to +\infty$,  this yields 
\begin{equation}\label{eq:4.3}
	\Big( \sum_{k=1}^\infty \Big\| \sum_{\lambda\in\Lambda_k} |e_\lambda^*(f)|  e_\lambda \Big\|_p^2  \Big)^{1/2} \le KC_p \|f\|_p.
\end{equation}

Next, similarly to \cite[Section 3.1]{LT25}, we claim that
\begin{equation}\label{eq:4.4}
	 \Big\| \sum_{\lambda\in\Lambda_k} |e_\lambda^*(f)|  e_\lambda \Big\|_p \ge \frac{1}{2}\sum_{\lambda\in\Lambda_k} |e_\lambda^*(f)|.
\end{equation}
Indeed, fix an arbitrary point $\mu_k \in Q_k$. 
Suppose that $\{a_\lambda\}$ are nonnegative
scalars  such that only finitely many of them are nonzero,
then 
\begin{align}
& \Big\| \sum_{\lambda\in\Lambda_k} a_\lambda e_\lambda \Big\|_p \ge  \Big\| \sum_{\lambda\in\Lambda_k} a_\lambda e_{\mu_k} \Big\|_p -  \Big\| \sum_{\lambda\in\Lambda_k} a_\lambda (e_{\mu_k} - e_\lambda) \Big\|_p \label{eq:LT.1} \\ 
 & \qquad =  \Big( \sum_{\lambda\in\Lambda_k} a_\lambda \Big) -   \Big\| \sum_{\lambda\in\Lambda_k} a_\lambda (e_{\mu_k} - e_\lambda) \Big\|_p
 \ge \frac{1}{2} \sum_{\lambda\in\Lambda_k} a_\lambda, \label{eq:LT.2}
\end{align}
where in \eqref{eq:LT.2} we used the fact that
$\|e_{\mu_k}\|_p = 1$ (since $\Om$ has measure $1$) and
that
 $|\mu_k - \lam| \le \delta$ and hence 
 $\|e_{\mu_k} - e_\lambda\|_p  < 1/2$.
If $\Lam_k$ happens to be a finite set, 
 we may apply \eqref{eq:LT.1}, \eqref{eq:LT.2}
 with $a_\lambda = |e_\lambda^*(f)|$ and arrive at \eqref{eq:4.4}.
In the  case where $\Lam_k$ is an infinite set, one can
establish \eqref{eq:4.4} by a  limiting argument,
 using the fact that the series
$\sum_{\lambda\in\Lambda_k}  |e_\lambda^*(f)| e_\lambda$
converges unconditionally.

Combining \eqref{eq:4.3}, \eqref{eq:4.4} we obtain the key estimate
\begin{equation}
\label{eq:key_estem}
	\Big( \sum_{k=1}^\infty \Big( \sum_{\lambda\in \Lambda_k} |e_{\lambda}^*(f)| \Big)^2 \Big)^{1/2} \le 2KC_p\|f\|_p
\end{equation}
which holds for every $f\in L^p(\Omega)$.

\subsection{}
We will now show that the estimate \eqref{eq:key_estem} implies that 
the series \eqref{eq:frame_expansion} in fact converges unconditionally
in $L^2(\Omega)$. 
Since there exist functions $f \in L^p(\Om) \setminus L^2(\Om)$, this
 will lead to the desired contradiction.

We begin by partitioning the family of cubes
 $\{Q_k\}_{k=1}^\infty$ into a finite number of
 sub\-families $\{Q_k\}$, $k\in A_j$, where $j=1,2,\dots,n$,
 in such a way  that
 \begin{equation}
 \label{eq:S.P.1}
 \dist(Q_k, Q_l) > 5\del, \quad k,l \in A_j, \quad k \ne l,
\end{equation}
that is,  any two cubes of the
 same subfamily are separated by distance at least $5 \del$.

Next, we fix a smooth function $\pphi$ on $\R^d$, supported in the ball $B_\del$
of radius $\del$ centered at the origin, such that its  Fourier transform
\begin{equation}
\ft{\varphi}(x)=\int_{\R^d} \varphi(t) e^{-2 \pi i \dotprod{x}{t}} dt
\end{equation}
satisfies 
  $|\ft{\pphi}(-x)| \ge 1$ for  $x\in\Omega$.
Such a function  $\pphi$
 exists since $\Om$ is a bounded set. 
 
 Now suppose that we fix $j\in\{1, 2,\ldots, n\}$ 
  and let $\{a_\lambda\}$ be
  a sequence of
scalars  such that only finitely many of them are nonzero.
Then we have
\begin{align}
&  \Big\| \sum_{k\in A_j}
  \sum_{\lambda\in\Lambda_k} a_\lambda  e_\lambda\Big\|_{L^2(\Om)}^2
 \le \int_{\R^d} \Big| \ft{\pphi}(-x)\sum_{k\in A_j}\sum_{\lambda\in\Lambda_k} a_{\lambda}  e_\lambda(x) \Big|^2 dx\\
& \qquad \qquad  =  \int_{\R^d}  \Big| \sum_{k\in A_j}\sum_{\lambda\in\Lambda_k} a_\lambda \pphi (t-\lambda) \Big|^2 dt, \label{eq:P.D.4.9}
\end{align}
where the last equality holds due to Plancherel's theorem.
We now observe that
since $\pphi$ is supported in $B_\del$ and due to \eqref{eq:S.P.1},
the functions $\sum_{\lambda\in\Lambda_k} a_\lam \pphi (t-\lambda)$,
$k \in A_j$, have disjoint supports. Therefore, the integral in \eqref{eq:P.D.4.9} is equal to
\begin{equation}
\sum_{k\in A_j} \int_{\R^d}
 \Big|\sum_{\lambda\in\Lambda_k} a_\lambda  \pphi(t-\lambda) \Big|^2 dt
 \le 
\|\pphi\|_2^2  \sum_{k\in A_j}  \Big(\sum_{\lambda\in\Lambda_k} |a_\lambda|  \Big)^2.
\end{equation}

To summarize, we have shown that
for  any $j \in \{1,2,\dots,n\}$, and for any sequence of scalars
$\{a_\lambda\}$ with only finitely many nonzero terms, we have
\begin{equation}
\label{eq:E.5.7}
  \Big\| \sum_{k\in A_j}
  \sum_{\lambda\in\Lambda_k} a_\lambda  e_\lambda\Big\|_{L^2(\Om)}
 \le 
M \Big(  \sum_{k\in A_j}  \Big(\sum_{\lambda\in\Lambda_k} |a_\lambda|  \Big)^2
\Big)^{1/2},
\end{equation}
where $M = \|\pphi\|_2$ is a constant not depending on 
the scalars $\{a_\lambda\}$.

Since the series \eqref{eq:frame_expansion} may be represented as
\begin{equation}
\label{eq:S.E.2}
	f=\sum_{j=1}^n \sum_{k\in A_j}  \sum_{\lambda\in\Lambda_k} e_\lambda^*(f) e_\lambda,
\end{equation}
then by applying \eqref{eq:E.5.7} for $j=1, 2, \ldots, n$ and 
using the estimate \eqref{eq:key_estem}, it  is easy to verify 
(for instance, using \cite[Proposition 2.4.9]{AK16})
that the series \eqref{eq:frame_expansion} in fact
converges unconditionally in $L^2(\Omega)$,
and that its sum, which ought to be the same function $f$,
satisfies the estimate $\|f\|_2 \le 2nMKC_p \|f\|_p$.
Since $p < 2$ this gives us a contradiction,
and the first part of
 \thmref{thm:noframes_exp} is therefore proved.
\qed


\section{Unconditional Schauder frames of exponentials: the case $p > 2$}

In this last section, we prove part \ref{nsf:ii} of \thmref{thm:noframes_exp}.

\subsection{}
Suppose that the set $\Omega\subset\R^d$ has an interior point. Then $\Omega$ contains a ball $B$. If $\{(e_\lambda, e_\lambda^*)\}$, $\lambda\in\Lambda$,
 is an unconditional Schauder frame for $L^p(\Omega)$, then the same system also forms an unconditional Schauder frame in $L^p(B)$. Hence we may assume that $\Omega$ itself is a ball.
Moreover, after translation and rescaling
we may assume that $\Om$ is
 a ball of measure $1$ centered at the origin.
 So we need to prove the following:

\begin{thm}
Let $\Omega = B_a\subset\R^d$
be a ball centered at the origin, $\mes(\Omega) = 1$. 
Then in the space $L^p(\Omega)$,  $p>2$, 
there is no unconditional Schauder frame of the form
 $\{(e_\lambda, e_\lambda^*)\}$, $\lambda\in\Lambda$,
  where $\Lambda\subset\R^d$ is a countable set and $\{e_\lambda^*\}$ are arbitrary elements of $(L^p(\Omega))^*$. 
\end{thm}

We recall that $B_a$ denotes the ball 
of radius $a$ centered at the origin.

\subsection{} 
As before, let $q=p/(p-1)$. Then $1<q<2$, and due to 
\lemref{lem:changeroles} the system
$\{(e^*_\lam, e_\lam)\}$,  $\lambda\in\Lambda$, forms
an unconditional Schauder frame in the space $L^q(\Om)$.
Since the action of the exponential
$e_\lam$, as an element of 
the space $L^p(\Om) = (L^q(\Om))^*$,
on a function  $g \in L^q(\Om)$  is given by 
$\ft{g}(- \lam)$, then every $g \in L^q(\Om)$ has an expansion
\begin{equation}
\label{eq:G.Q.1}
g = \sum_{\lam\in\Lambda} \ft{g}(-\lam) e^*_\lam
\end{equation}
unconditionally convergent in $L^q(\Om)$. 
By \lemref{lem:unc_constant}, there is a constant $K$ such that
\begin{equation}
\label{eq12}
\Big\| \sum_{\lam \in \Lam} \theta_\lam \ft{g}(-\lam) e^*_\lam \Big\|_q \le K  \|g\|_q
\end{equation}
holds for any sequence  of scalars $\{\theta_\lam\}$ with 
$|\theta_\lam| \le 1$.

\subsection{}
Let $Q = [-\tfrac1{2}, \tfrac1{2})^d$
be the half-open unit cube centered at the origin,
and  consider a partition of $\R^d$ 
into cubes $\{ Q + k \}$, $k \in \Z^d$.
It induces a 
partition of $\Lam$ into blocks  $\Lam_k = \Lam \cap (Q+k)$,
$k \in \Z^d$. Again we keep in mind that each block $\Lam_k$
may contain arbitrarily many
elements, or  may  even be infinite.

To each $g \in L^q(\Om)$ we associate a sequence
$\{w_k(g)\}$, $k \in \Z^d$, defined by
\begin{equation}
\label{eq32}
w_k(g) = \sup 
\Big\| \sum_{\lam \in \Lam_k} \theta_\lam \ft{g}(-\lam) e^*_\lam \Big\|_q
\end{equation}
where the supremum is taken over all sequences
of scalars $\{\theta_\lam\}$ with 
$|\theta_\lam| \le 1$.  We claim that the estimate
\begin{equation}
\label{eq6}
\Big( \sum_{k \in \Z^d} w_k(g)^2 \Big)^{1/2} \le  K C_q  \|g\|_q
\end{equation}
holds for every $g \in L^q(\Om)$. 

Indeed, let $\eps > 0$, and 
for each $k\in\Z^d$ choose the scalars 
$\{\theta_\lambda\}$, $\lambda\in\Lambda_k$, so that
\begin{equation}
\label{eq:5.3}
w_k(g)\le (1+\eps)
\Big\| \sum_{\lam \in \Lam_k} \theta_\lam \ft{g}(-\lam) e^*_\lam \Big\|_q.
\end{equation}
By \eqref{eq12}, for every $N$ and every $t\in [0,1]$ we have 
\begin{equation}
\Big\| \sum_{|k| \le N} r_{n_k}(t) \sum_{\lambda\in\Lambda_k} \theta_\lam \ft{g}(-\lam) e^*_\lam \Big\|_q\le K\|g\|_q,
\end{equation}
where $\{r_{n_k}(t)\}$, $k \in \Z^d$, is some  enumeration of the
 Rademacher functions. 
If we now take $\int_0^1 dt$, apply \lemref{lem:cotype}, use \eqref{eq:5.3},
let $N \to + \infty$ and  $\eps \to 0$, we arrive at \eqref{eq6}.

\subsection{} 
For any function $g\in L^q(\R^d)$ supported in $\Omega = B_a$, we have
\begin{equation}
	\label{eq1}
	\|\ft{g}\|_\infty \le \|g\|_1 \le \|g\|_q 
	\le  \sum_{k \in \Z^d} \Big\| \sum_{\lam \in \Lam_k}
	\ft{g}(-\lam)   e^*_\lam \Big\|_{q},
\end{equation}
where the second inequality holds since $\mes(\Om)=1$,
while in the last inequality we used the
unconditionally convergent series
expansion \eqref{eq:G.Q.1} of the function $g$.

\subsection{} 
Fix any $0 < r < a$. 
We claim that there is a constant $M=M(d,a,r)$ such that
\begin{equation}
\label{eq2}
\|\ft{f}\|_2 \leq M \Big(
\sum_{k \in \Z^d} w_k(f)^2 \Big)^{1/2}
\end{equation}
holds for every 
$f \in L^2(\R^d)$ with $\supp(f) \sbt B_r$.
Indeed, fix a smooth function $\varphi$ 
such that
$\supp(\varphi) \sbt B_{a-r}$ and 
$\ft{\varphi}(0) = 1$. 
Define $h_x(y) = \ft{f}(y) \ft{\varphi}(x-y)$, then
it is easy to see that
$h_x = \ft{g}_x$ where $g_x$ is 
a smooth function  supported in $\Om = B_a$.
In particular the estimate \eqref{eq1} can be
applied to  $g = g_x$, which yields
\begin{equation}
	\label{eq3}
	|\ft{f}(x)| = |h_x(x)|  
	\le  \sum_{k \in \Z^d}\Big\| \sum_{\lam \in \Lam_k}
	h_x(-\lam)   e^*_\lam \Big\|_{q}
	= \sum_{k \in \Z^d} \Big\| \sum_{\lam \in \Lam_k}
	\ft{f}(-\lam) \ft{\varphi}(x + \lam)   e^*_\lam \Big\|_{q}.
\end{equation}
Now set
\begin{equation}
	\label{eq38}
	\Phi(x) = \sup_{y \in Q} |\ft{\varphi}(x+y)|, \quad x \in \R^d,
\end{equation}
and note that
$| \ft{\varphi}(x + \lam) | \le \Phi(x+k)$ for $\lam \in \Lam_k$.
Hence,
using \eqref{eq32}, \eqref{eq3}  we obtain
\begin{equation}
	\label{eq31}
	|\ft{f}(x)| \le \sum_{k \in \Z^d} w_k(f) \cdot \Phi(x+k).
\end{equation}
An application of the Cauchy-Schwarz inequality now yields
\begin{equation}
	\label{eq4}
	|\ft{f}(x)|^2 \le   \Big( \sum_{k \in \Z^d} w_k(f)^2 \cdot \Phi(x+k) \Big)
	\Big( \sum_{k \in \Z^d} \Phi(x+k) \Big).
\end{equation}
Since the function $\Phi$ has fast decay,
both $\sup_x \sum_{k}  \Phi(x + k)$ and
$\int \Phi(x)dx$ do not exceed a certain constant $M$.
Hence taking $\int_{\R^d} dx$ of the inequality \eqref{eq4} yields \eqref{eq2}.

\subsection{}
Finally let $f \in L^2(\R^d)$,  $\supp(f) \sbt B_r$.
Combining \eqref{eq6} and \eqref{eq2}   yields
\[
\|f\|_2 = \|\ft{f}\|_2 
 \le M \Big( \sum_{k \in \Z^d}
w_k(f)^2 \Big)^{1/2} \le M K C_q \|f\|_q,
\]
where the first equality holds due to Plancherel's theorem.
Since $q < 2$, the latter inequality leads to the desired contradiction,
and so the second part of
 \thmref{thm:noframes_exp} is now proved.
\qed


\end{document}